\documentclass[11pt]{article}
\usepackage{amsmath}
 \usepackage{a4}                          
\usepackage{amssymb}
\DeclareMathOperator{\Log}{Log}

\newcommand{\Int}{\int_{-\infty}^\infty}
 \newcommand{\eps}{\varepsilon}
 
 \renewcommand{\a}{\alpha}
 \newcommand{\La}{\Lambda}

\newtheorem{thm}{Theorem}[section]
\newtheorem{lemma}[thm]{Lemma}

\newtheorem{rem}[thm]{Remark}
\newtheorem{cor}[thm]{Corollary}

\title{Asymptotic behaviour of the Urbanik semigroup}

\author{Christian Berg\thanks{The  first author has
    been supported by grant 10-083122 from  {\it The Danish
    Council for Independent Research $|$ Natural Sciences}}
\and Jos\'e Luis L\'opez
 \thanks{The second author has been supported by  grant MTM2010-21037 from the {\it Direcci\'on General de Ciencia y Tecnolog\'{\i}a}}
}
\begin{document}
\maketitle

\begin{abstract} We revisit the product convolution semigroup of probability densities $e_c(t),\,c>0$
on the positive half-line with moments $(n!)^c$ and determine the asymptotic behaviour of $e_c$ for large and small $t>0$. This shows that $(n!)^c$ is indeterminate as Stieltjes moment sequence if and only if $c>2$.
  \end{abstract}


\section{Introduction} We consider a family of
probability densities $e_c(t),c>0$ on the half-line given by 
\begin{equation}\label{eq:U1}
e_c(t)=\frac{1}{2\pi}\Int t^{ix-1}\Gamma(1-ix)^c\,dx,\quad t>0.
\end{equation}
In this formula we use that $\Gamma(z)$ is a non-vanishing holomorphic
function in the cut plane 
\begin{equation}\label{eq:cut}
\mathcal A=\mathbb C\setminus (-\infty,0],
\end{equation}
so we can define
$$
\Gamma(z)^c=\exp(c\log\Gamma(z)),\quad z\in \mathcal A
$$
using the holomorphic branch of $\log\Gamma$ which is 0 for $z=1$.

As far as we know it was proved first by Urbanik in \cite[Section 4]{U1} that
$e_c$ is a probability density, and that the following product
convolution equation holds
\begin{equation}\label{eq:U2}
  e_{c+d}(t)=\int_0^\infty e_c(t/x)e_d(x)\frac{dx}{x},\quad c,d>0.
\end{equation}
Furthermore, it was noticed that
\begin{equation}\label{eq:U3}
\int_0^\infty t^ne_c(t)\,dt=(n!)^c,\quad c>0,n=0,1,\ldots.
\end{equation}
Defining the probability measure $\tau_c$ on $(0,\infty)$ by
\begin{equation}\label{eq:mea1}
d\tau_c=e_c(t)\,dt=te_c(t)\,dm(t),\quad c>0,
\end{equation}
where $dm(t)=(1/t)\,dt$ is the Haar measure on the locally compact
abelian group $G=(0,\infty)$ under multiplication, we can write
\eqref{eq:U2} as $\tau_c\diamond \tau_d=\tau_{c+d}$, where
$\diamond$ denotes the (product) convolution of measures on the
multiplicative group $G$. The family  
$(\tau_c)_{c>0}$ is a convolution semigroup in the sense of \cite{B:F}.
We propose to call this semigroup the Urbanik semigroup because of
\cite{U1}.

The continuous characters of the group $G$  can be given as $t\to
t^{ix}$, where $x\in\mathbb R$ is arbitrary, and in this way the dual
group $\widehat{G}$ of $G$ can be identified with the additive group
of real numbers, and by the inversion theorem of Fourier analysis for
LCA-groups, \eqref{eq:U1} is equivalent to 
\begin{equation}\label{eq:U4}
\widehat{\tau_c}(x)=\int_0^\infty t^{-ix}\,d\tau_c(x)=\exp(c\log(\Gamma(1-ix)),\quad x\in\mathbb R.
\end{equation}
To establish the existence of a product convolution semigroup
$(\tau_c)$ satisfying \eqref{eq:U4} is therefore equivalent to
proving that
\begin{equation}\label{eq:U5}
\rho(x):=-\log\Gamma(1-ix),\quad x\in \mathbb R
\end{equation}
 is a continuous negative definite function on $\mathbb R$ in the
 terminology of \cite{B:F} or \cite{S:S:V}.

This was done in \cite{U1} by giving the L\'evy-Khinchin
representation of $\rho$, using Malmsten's formula, cf. \cite[8.341(3)]{G:R}:
\begin{equation}\label{eq:Malmsten}
\log\Gamma(z)=\int_0^\infty\left[\frac{e^{-zt}-e^{-t}}{1-e^{-t}}+(z-1)e^{-t}\right]\,\frac{dt}{t},\quad \textrm{Re\,}(z)>0.
\end{equation}
 
In fact this formula can be written
\begin{equation}\label{eq:U6}
-\log\Gamma(1-ix)=\int_0^\infty\left[1-e^{ixt}+\frac{itx}{1+t^2}\right]\frac{e^{-t}}{t(1-e^{-t})}\,dt-iax,
\end{equation}
where
$$
a=\int_0^\infty\left[\frac{1}{(1+t^2)(1-e^{-t})}-\frac{1}{t}\right]e^{-t}\,dt,
$$
showing that $\rho(x)=-\log\Gamma(1-ix)$ is negative definite with the
L\'evy measure 
$$
d\mu=\frac{e^{-t}}{t(1-e^{-t})}\,dt
$$
concentrated on $(0,\infty)$.

Another proof of the negative definiteness of $\rho$ was given in \cite{B} based
on the Weierstrass product for $\Gamma$ , where $\Log$ denotes the principal logarithm in the cut plane $\mathcal A$, cf. \eqref{eq:cut}:
$$
-\log\Gamma(z)=\gamma z + \Log z
+\sum_{k=1}^\infty\left(\Log(1+z/k)-z/k\right),\quad z\in\mathcal A. 
$$
Clearly,
$$
\rho_n(z):=\gamma z + \Log z
+\sum_{k=1}^n\left(\Log(1+z/k)-z/k\right) 
$$
converges locally uniformly to $-\log\Gamma(z)$ for $z\in\mathcal A$,
and since
$$
\rho_n(1-ix)=\rho_n(1) -i\left(\gamma-\sum_{k=1}^n\frac1k\right)x +\sum_{k=1}^{n+1}\Log(1-ix/k)
$$
is negative definite, because  $\Log(1+iax)$ is so for $a\in\mathbb R$ and
$$
\rho_n(1)=\gamma+\log(n+1)-\sum_{k=1}^n\frac1k > 0,
$$
we conclude that the limit function $\rho(x)=-\log\Gamma(1-ix)$ is
negative definite.

As noticed in \cite[Lemma 2.1]{B}, \eqref{eq:U3} is a special case of
\begin{equation}\label{eq:hol}
\int_0^\infty t^ze_c(t)\,dt=\Gamma(1+z)^c,\quad \textrm{Re\,}(z)>-1,
\end{equation} 
and letting $z$ tend to $-1$ along the real axis, we get
\begin{equation}\label{eq:infmass}
\int_0^\infty e_c(t)\frac{dt}{t}=\int_0^\infty e_c(1/t)\frac{dt}{t}=\infty,\quad c>0.
\end{equation}

It follows from \eqref{eq:U3} that $(n!)^c$ is a Stieltjes moment
sequence for any $c>0$, and while it is easy to see that it is
S-determinate for $c\le 2$ in the sense, that there is only one measure
on the half-line with these moments, namely $\tau_c$, it is rather
delicate to see that it is S-indeterminate for $c>2$. This was proved in
Theorem 2.5 in \cite{B}. The proof was based on a relationship between
$\tau_c$ and stable distributions, and it used heavily asymptotic
results of Skorokhod from \cite{Sk} and exposed in \cite{Z}. Further details are given at the
end of this section.

The purpose of the present paper is to establish the
asymptotic behaviour of the densities $e_c(t)$ for $t\to \infty$ and $t\to 0$.
The behaviour for $t\to\infty$ will lead to a direct proof of the S-indeterminacy for $c>2$.
 
We mention that the product convolution semigroup $(\tau_c)_{c>0}$ corresponds to the Bernstein function $f(s)=s$ in the following result from \cite[Theorem 1.8]{B}.

\begin{thm}\label{thm:infdiv} Let $f$ be a non-zero Bernstein function.
The uniquely determined measure $\kappa=\kappa(f)$ with moments
$s_n=f(1)\cdots f(n)$  is infinitely divisible with respect to
the product convolution.  The unique product convolution semigroup
 $(\kappa_c)_{c>0}$  with $\kappa_1=\kappa$ has the moments
\begin{equation}\label{eq:pcsmoments}
 \int_0^\infty x^n\,d\kappa_c(x)=(f(1)\cdots f(n))^c,\quad c>0, n=0,1,\ldots.
\end{equation}
\end{thm}

It is an easy consequence of Carleman's criterion that  the measures $\kappa_c$ are S-determinate
for $c\le 2$, cf. \cite[Theorem 1.6]{B}.
  
In \cite{B} we consider three Bernstein functions $f_\alpha,f_\beta,
f_\gamma$ with corresponding product convolution semigroups $(\alpha_c)_{c>0},(\beta_c)_{c>0},
(\gamma_c)_{c>0}$:
$$
f_\a(s)=(1+1/s)^s,\quad f_\beta(s)=(1+1/s)^{-s-1},\quad f_\gamma(s)=s(1+1/s)^{s+1}.
$$
It is proved that  the measures $\alpha_c,\beta_c$ have compact support, so they are clearly S-determinate for all $c>0$, but $\gamma_c$ is S-indeterminate for $c>2$. Using that
$\tau_c=\beta_c\diamond\gamma_c$, it is possible to infer that also $\tau_c$ is S-indeterminate,
see \cite{B} for details.

As noticed in \cite{U1}, the measures $\tau_c,\;c\ge 1$ are also infinitely divisible for 
the additive structure, because $e_c(t)$ is completely monotonic. To see this, notice that 
the convolution equation \eqref{eq:U2} with $d=1$ can be written
\begin{equation}\label{eq:cm}
e_{c+1}(t)=\int_0^\infty e^{-tx}e_c(1/x)\,\frac{dx}{x},\quad c>0,
\end{equation}
showing that $e_c(t)$ is completely monotonic for $c>1$, and it tends to infinity for $t\to 0$ because of \eqref{eq:infmass}.  

It is well-known that the exponential distribution $\tau_1$ is infinitely divisible for the additive structure and with a completely monotonic density $e_1(t)$.

 Urbanik also showed that $\tau_c$ is not infinitely divisible for the additive structure when $0<c<1$.

Formula \eqref{eq:U1} states roughly speaking that $te_c(t)$ is the Fourier transform  of the Schwartz function $\Gamma(1-ix)^c$ evaluated at $\log t$, thus showing that $e_c$
is $C^\infty$ on $(0,\infty)$. By Riemann-Lebesgue's Lemma we also see that 
$te_c(t)$ tends to zero for $t$ tending to zero and to infinity. Much more will be obtained in the main results below.

\section{Main results}
Our main results are

\begin{thm}\label{thm:main} For $c>0$ we have
\begin{equation}\label{eq:as1}
e_c(t) = \frac{(2\pi)^{(c-1)/2}}{\sqrt{c}}\frac{\exp(-ct^{1/c})}{t^{(c-1)/(2c)}}
\left[1+\mathcal O\left(\frac{1}{t^{1/c}}\right)\right],
\quad t\to\infty. 
\end{equation}
\end{thm}

\begin{rem}\label{thm:remark1} {\rm The densities $e_c$ are not explicitly known except for $c=1,2$, where 
$$
e_1(t)=e^{-t},\quad e_2(t)=\int_0^\infty\exp(-x-t/x)\frac{dx}{x}=2K_0(2\sqrt{t}).
$$
In the last formula $K_0$ is a modified Bessel function, see \cite[Chap. 10, Sec. 25]{O:M}.}
\end{rem}

\begin{cor}\label{thm:indet}  
The measure $\tau_c=e_c(t)\,dt$ is S-indeterminate for $c>2$.
\end{cor}

\begin{thm}\label{thm:main2} For $c>0$ we have
\begin{equation}\label{eq:asat0}
e_c(t)=\frac{(\log(1/t))^{c-1}}{\Gamma(c)} +\mathcal O((\log(1/t))^{c-2}),\quad t\to 0.
\end{equation}
\end{thm}

\begin{rem}\label{thm:remark2} {\rm Formula \eqref{eq:asat0} shows that $e_c(t)$ tends to infinity as a power of $\log(1/t)$ when $c>1$, but so slowly that multiplication with $t$ forces the density to tend to zero. When $0<c<1$ the density $e_c(t)$ tends to zero.}
\end{rem}

\section{Proofs} We will first give a proof of Theorem~\ref{thm:main} in the case, where $c$ is a natural number. Note that the asymptotic expression in \eqref{eq:as1} for $c=1$ reduces to
$e_1(t)=e^{-t}$. When $c=n+1$, where $n$ is a natural number, we know that $e_{n+1}(t)$ is the $n$'th product convolution power of $e_1$, hence
$$
e_{n+1}(t)=\int_0^\infty\ldots\int_0^\infty\,e^{-\frac{t}{u_1\cdots u_n}}e^{-u_1}\cdots e^{-u_n}
\,\frac{du_1}{u_1}\cdots\frac{du_n}{u_n}.
$$
For $t>0$  fixed,  the change of variables $u_j=t^{1/(n+1)}v_j,j=1,\ldots,n$ leads to
\begin{equation}\label{eq:as2}
e_{n+1}(t)=\int_0^\infty \ldots\int_0^\infty g(v_1,\ldots,v_n)e^{-t^{1/(n+1)}f(v_1,\ldots,v_n)}dv_1\cdots dv_n,
\end{equation}
with
$$
g(v_1,\ldots,v_n):=\frac{1}{v_1\cdots v_n},\quad f(v_1,\ldots,v_n):=v_1+\cdots+v_n+g(v_1,\ldots,v_n).
$$
The phase function $f(v_1,\ldots,v_n)$ is convex in $C=\{v_1>0,\ldots,v_n>0\}$ because the
Hessian matrix of second derivatives  is
$$
Hf(v_1,\ldots,v_n)=g(v_1,\ldots,v_n)\left(\begin{matrix}
\frac{2}{v_1^2} & \frac{1}{v_1v_2} &\cdots &\frac{1}{v_1v_n}\\
\frac{1}{v_2v_1} & \frac{2}{v_2^2} &\cdots &\frac{1}{v_2v_n}\\
\vdots &\vdots & \ddots  & \vdots\\
\frac{1}{v_nv_1} & \frac{1}{v_nv_2} & \cdots &\frac{2}{v_n^2}
\end{matrix}\right),
$$
which is easily seen to be positive definite. The phase function therefore has a global minimum at the unique stationary point  $\vec{v}_0$ such that $\vec{\bigtriangledown}f(\vec{v}_0)=\vec{0}$, that is, at $\vec{v}_0=(1,\ldots,1)$.  At that point, the Hessian matrix of $f(\vec{v})$ is
$$
A:=Hf(1,\ldots,1)=\left(\begin{matrix} 2 & 1 & 1& \cdots & 1 \\ 1 & 2 & 1 & \cdots
 & 1 \\
\vdots &\vdots & \vdots  & \ddots & \vdots \\
1 & 1 & 1 & \cdots & 2 
\end{matrix}\right),
$$
with determinant $\det(A)=n+1$.

By Laplace's asymptotic method for multiple dimensional Laplace transforms, cf. \cite[Theorem 3, p. 495]{W}, we know that for $t\to\infty$,
$$
e_{n+1}(t) = \left(\frac{2\pi}{t^{1/(n+1)}}\right)^{n/2} g(\vec{v}_0)(\det(A))^{-1/2} e^{-t^{1/(n+1)}f(\vec{v}_0)}\left[1+{\mathcal O}\left(\frac{1}{t^{1/(n+1)}}\right)\right].
$$
We have that $g(\vec{v}_0)=1$ and $f(\vec{v}_0)=n+1$, hence

\begin{equation}
e_{n+1}(t)  =
\frac{(2\pi)^{n/2}}{\sqrt{n+1}}\frac{e^{-(n+1)t^{1/(n+1)}}}{t^{n/(2(n+1))}}
\left[1+{\mathcal O}\left(\frac{1}{t^{1/(n+1)}}\right)\right],
\end{equation}
which agrees with \eqref{eq:as1} for $c=n+1$.

\medskip
The proof of Theorem~\ref{thm:main} for arbitrary $c>0$ is more delicate. We first apply Cauchy's integral theorem to move the integration in \eqref{eq:U1} to an arbitrary horizontal line
\begin{equation}\label{eq:line}
L_a:=\{z=x+ia\mid x\in\mathbb R\},\quad a>0.
\end{equation}

\begin{lemma}\label{thm:shift} With $L_a$ as in \eqref{eq:line} we have
\begin{equation}\label{eq:def}
e_c(t)=\frac{1}{2\pi}\int_{L_a} t^{iz-1}\Gamma(1-iz)^c\,dz,\quad t>0.
\end{equation}
\end{lemma}

{\it Proof:} For $t,c>0$ fixed, $f(z)=t^{iz-1}\Gamma(1-iz)^c$ is holomorphic in the simply connected domain $\mathbb C\setminus i(-\infty,-1]$, so the Lemma follows from Cauchy's integral theorem provided the integral
$$
\int_0^{a} f(x+iy)\,dy
$$
tends to 0 for $x\to\pm\infty$. We have
$$
|f(x+iy)|=t^{-y-1}|\Gamma(1+y-ix)|^c
$$
and since
$$
 |\Gamma(u+iv)|\sim \sqrt{2\pi}e^{-\pi/2|v|}|v|^{u-1/2},\quad |v|\to\infty, \;\mbox{uniformly for bounded real $u$},
$$
cf. \cite[p.141, eq. 5.11.9]{NISTA:R},\cite[8.328(1)]{G:R}, the result follows. $\quad\square$

\medskip
In the following we will use  Lemma~\ref{thm:shift} with the line of integration $L=L_a$, where $a=t^{1/c}-1$ for $t>1$.
Therefore, using the parametrization $z=x+i(t^{1/c}-1)$ we get
$$
e_c(t)=t^{-t^{1/c}}\frac{1}{2\pi}\Int t^{ix}\Gamma(t^{1/c}-ix)^c\,dx,
$$
and after the change of variable $x=t^{1/c}u$
\begin{equation}\label{eq:def2}
e_c(t)=t^{1/c-t^{1/c}}\frac{1}{2\pi}\Int t^{iut^{1/c}}\Gamma(t^{1/c}(1-iu))^c\,du.
\end{equation}

Binet's formula for $\Gamma$ is (\cite[8.341(1)]{G:R})
\begin{equation}\label{eq:Binet}
\Gamma(z)=\sqrt{2\pi}z^{z-\frac12}e^{-z+\mu(z)},\quad \textrm{Re\,}(z)>0,
\end{equation}
where
\begin{equation}\label{eq:Binet1}
\mu(z)=\int_0^\infty\left(\frac{1}{2}-\frac{1}{t}+\frac{1}{e^t-1}\right)\frac{e^{-zt}}{t}\,dt,\quad \textrm{Re\,}(z)>0.
\end{equation}
Notice that $\mu(z)$ is the Laplace transform of a positive function, so we
have the estimates for $z=r+is, r>0$
\begin{equation}\label{eq:Binet2}
|\mu(z)|\le \mu(r)\le \frac{1}{12r},
\end{equation}
where the last inequality is a classical version of Stirling's formula, thus showing that the estimate is uniform in $s\in\mathbb R$.

Inserting this in \eqref{eq:def2}, we get after some simplification
\begin{equation}\label{eq:def3}
e_c(t)=(2\pi)^{c/2-1}t^{1/c-1/2}e^{-ct^{1/c}}\Int e^{ct^{1/c}f(u)}g_c(u)M(u,t)\,du,
\end{equation}
where
\begin{equation}\label{eq:hp1}
f(u):=iu+(1-iu)\Log(1-iu),\quad g_c(u):=(1-iu)^{-c/2}
\end{equation}
and
\begin{equation}\label{eq:hp2}
M(u,t):=\exp[c\mu(t^{1/c}(1-iu))].
\end{equation}
From \eqref{eq:Binet2} we get $M(u,t)=1+\mathcal O(t^{-1/c})$ for
$t\to\infty$, uniformly in $u$.
We shall therefore consider the behaviour of
\begin{equation}\label{eq:final}
\Int  e^{ct^{1/c}f(u)}g_c(u)\,du.
\end{equation}
From here we need to apply the saddle point method to obtain the
 approximation of \eqref{eq:final} for large positive
$t$. For convenience, we use Theorem 1 in \cite{L:P:P}. We have that
the only saddle point of the phase function $f(u)$ is $u=0$ and
$f(0)=f'(0)=0$, $f''(0)=-1$, $f'''(0)\ne 0$; also $g_c(0)=1$. Then,
the parameters used in that theorem are $m=2$, $p=3$, $\phi=\pi$,
$N=0$, $M=1$ and the large variable used in the theorem is $x\equiv
ct^{1/c}$. We have that the steepest descendent path used in the
theorem is
$\Gamma=\Gamma_0\bigcup\Gamma_1=(-\infty,0)\bigcup(0,\infty)$, that
is, it is just the original integration path in the above integral,
and therefore does not need any deformation. From \cite[Theorem
1]{L:P:P} with the notation used there,
we read that the integral \eqref{eq:final} has an expansion of the form
$$
e^{xf(0)}[c_0\Psi_0(x)+c_1\Psi_1(x)+c_2\Psi_2(x)+\cdots],
$$
with $\Psi_{n}(x)=\mathcal O(x^{-(n+1)/2})$ and $c_n$ is independent of $x$.
Because the factors $c_{2n+1}$ vanish we find
$$
c_0\Psi_0(x)+c_1\Psi_1(x)+c_2\Psi_2(x)+\cdots=c_0\Psi_0(x)[1+\mathcal O(x^{-1})]
$$
with $c_0=1$ and
$$
\Psi_0(x)=a_0(x)\Gamma\left(\frac12\right)\left\vert\frac{2}{x
    f''(0)}\right\vert^{1/2}
$$
with
$$
a_0(x)=e^{-xf(0)}A_0(x)B_0, \quad A_0(x)=e^{xf(0)}, \quad B_0=g_c(0),
$$
hence $a_0(x)=B_0=1$.
Using all these data we finally obtain
$$
\int_{-\infty}^\infty e^{ct^{1/c}f(u)}g_c(u)du=\frac{\sqrt{2\pi}}{\sqrt{c}t^{1/(2c)}}[1+{\mathcal O}(t^{-1/c})],
$$
and
$$
e_c(t)=\frac{(2\pi)^{(c-1)/2}}{\sqrt{c}}\frac{e^{-ct^{1/c}}}{t^{(c-1)/(2c)}}[1+{\mathcal O}(t^{-1/c})].
$$
$\square$

\medskip
{\it Proof of Corollary~\ref{thm:indet}.} We apply the Krein criterion for S-indeterminacy of probability densities concentrated on the half-line, using a version given in \cite[Theorem 5.1]{B0}. 
It states that if 
\begin{equation}\label{eq:as3}
\int_0^\infty \frac{\log e_c(t)\,dt}{\sqrt{t}(1+t)}>-\infty,
\end{equation}
then $\tau_c=e_c(t)\,dt$ is S-indeterminate. We shall see that
\eqref{eq:as3} holds for $c>2$.

From Theorem~\ref{thm:main} combined with the fact that $e_c(t)$ is decreasing when $c>1$,
we see that the inequality in \eqref{eq:as3} holds 
if and only if
$$
\int_0^\infty \frac{\log((2\pi)^{(c-1)/2}/\sqrt{c}) -ct^{1/c}-((c-1)/(2c))\log t}{\sqrt{t}(1+t)}\,dt>-\infty,
$$
and the latter holds precisely for $c>2$. This shows that $\tau_c$ is S-indeterminate for $c>2$.
$\quad\square$

\medskip
{\it Proof of Theorem~\ref{thm:main2}.}

Since we are studying the behaviour for $t\to 0$, we assume that
$0<t<1$ so that $\La:=\log(1/t)>0$.

We will need integration along the vertical lines
\begin{equation}\label{eq:vert1}
V_a:=\{a+iy\mid y=-\infty\ldots\infty\},\quad a\in \mathbb R,
\end{equation}
and we can therefore  express \eqref{eq:U1} as
\begin{equation}\label{eq:fund1}
e_c(t)=\frac{1}{2\pi i}\int_{V_{-1}} t^{z}\Gamma(-z)^cdz.
\end{equation}
By the functional equation for $\Gamma$ we get
\begin{equation}\label{eq:cont1}
e_c(t)=\frac{1}{2\pi i}\int_{V_{-1}} (-z)^{-c}t^{z}\Gamma(1-z)^cdz.
\end{equation}
To ease the writing we define
$$
\varphi(z):=t^z\Gamma(1-z)^{c},\quad g(z):=(-z)^{-c}=\exp(-c\Log(-z)),
$$
and note that $\varphi$ is holomorphic in $\mathbb C\setminus[1,\infty)$, while $g$
is holomorphic in $\mathbb C\setminus [0,\infty)$. Here $\Log$ is the
principal logarithm in the cut plane $\mathcal A$, cf. \eqref{eq:cut}. 

Note that for $x>0$ 
$$
g_{\pm}(x):=\lim_{\eps\to 0}g(x \pm i\eps)=x^{-c}e^{\pm i\pi c}.
$$

Formula \eqref{eq:cont1} can now be written
\begin{equation}\label{eq:cont2}
e_c(t)=\frac{1}{2\pi i}\int_{V_{-1}} g(z)\varphi(z)\,dz.
\end{equation}

\medskip
{\bf Case 1.} We will first treat the case $0<c<1$.

We fix $0<s<1$ and choose $0<\eps<s$  and integrate $g(z)\varphi(z)$ over the contour
$\mathcal C$
$$
\{-1+iy \mid y=\infty\ldots 0\}\cup [-1,-\eps]\cup \{\eps e^{i\theta}\mid\theta=\pi\ldots 0\}\cup[\eps,s]\cup\{s+iy\mid y=0\ldots \infty\}
$$ 
and get 0 by the integral theorem of Cauchy. On the interval $[\eps,s]$ we use the values of $g_+$.

Similarly we get 0 by integrating  $g(z)\varphi(z)$ over the complex conjugate contour $\overline{\mathcal C}$,
and now we use the values of $g_{-}$ on the interval $[\eps,s]$.

Subtracting the second contour integral from the first leads to
$$ 
\int_{V_s}-\int_{V_{-1}}-\int_{|z|=\eps} g(z)\varphi(z)\,dz + \int_\eps^s \varphi(x)
(g_{+}(x)-g_{-}(x))\,dx=0,
$$
where the integral over the circle is with positive orientation. Note
that the two integrals over $[-1,-\eps]$ cancel. Using that $0<c<1$ it is easy to see that the just mentioned integral converges to $0$ for $\eps\to 0$, and we finally get for $\eps\to 0$ 
 $$
e_{c}(t)=\frac{1}{2\pi i}\int_{V_s} g(z)\varphi(z)\,dz +
\frac{\sin(\pi c)}{\pi}\int_0^s  x^{-c}\varphi(x)\,dx:=I_1+I_2.
$$
We claim that the first integral $I_1$ is $o(t^s)$ for $t\to 0$. To
see this we insert the parametrization of $V_s$ and get
\begin{eqnarray*}  
I_1=\frac{t^s}{2\pi}\Int (-s-iy)^{-c}t^{iy}\Gamma(1-s-iy)^{c}\,dy
\end{eqnarray*}
and the integral is $o(1)$ by Riemann-Lebesgue's Lemma, so $I_1=o(t^s)$.

The substitution $u=x\log(1/t)=x\La$ in the integral $I_2$ leads to
\begin{equation}\label{eq:I2}  
I_2=\frac{\sin(\pi c)}{\pi}\La^{c-1}
\int_0^{s\La} u^{-c}e^{-u}\Gamma(1-u/\La)^c\,du.
\end{equation}
We split the integral in \eqref{eq:I2} as
\begin{equation}\label{eq:I2a}
\Gamma(1-c)+\int_0^{s\La} u^{-c}e^{-u}\left[\Gamma(1-u/\La)^c -1\right]\,du
-\int_{s\La}^\infty u^{-c}e^{-u}\,du,
\end{equation}
and by the mean-value theorem and $\Psi=\Gamma'/\Gamma$ we have 
$$
\Gamma(1-u/\La)^c -1=-\frac{u}{\La}c\Gamma(1-\theta
u/\La)^c\Psi(1-\theta u/\La)
$$
for some $0<\theta<1$, but this implies that
$$
|\Gamma(1-u/\La)^c-1|\le \frac{cu}{\La}M(s),\quad
0<u<s\La,
$$
where
$$
M(s):=\max\{\Gamma(x)^c|\Psi(x)|\mid 1-s\le x\le 1\},
$$
so the first integral in \eqref{eq:I2a} is $\mathcal O(\La^{-1})$.
The second integral is an incomplete Gamma function, and by known
asymptotics for this, see \cite{G:R}, we get that the second integral
is $\mathcal O(\La^{-c}t^s)$. Putting things together and
using Euler's reflection formula for $\Gamma$, we see
that
$$
e_c(t)=\frac{\La^{c-1}}{\Gamma(c)}+\mathcal
O(\La^{c-2}),
$$
which is \eqref{eq:asat0}.

\bigskip
{\bf Case 2.} We now assume $1<c<2$.

The Gamma function decays so rapidly when $z=-1+iy\in
V_{-1},y\to\pm\infty$, that we can integrate by parts in
\eqref{eq:cont1} to get
\begin{equation}\label{eq:cont3}
e_{c}(t)=-\frac{1}{2\pi i}\int_{V_{-1}} \frac{(-z)^{-(c-1)}}{c-1}\frac{d}{dz}(t^z\Gamma(1-z)^c)\,dz.
\end{equation}
Defining
$$
\varphi_1(z):=\frac{d}{dz}(t^z\Gamma(1-z)^c)=t^z\Gamma(1-z)^c(\log t-c\Psi(1-z)),
$$
and using the same contour technique as in case 1 to the integral in \eqref{eq:cont3}, where now $0<c-1<1$, we get for $0<s<1$ fixed
$$
e_c(t)=-\frac{1}{c-1}\frac{1}{2\pi i}\int_{V_s}(-z)^{-(c-1)}\varphi_1(z)\,dz -\frac{\sin(\pi(c-1))}{(c-1)\pi}\int_0^s x^{-(c-1)}\varphi_1(x)\,dx.
$$
The first integral is $o(t^s\La)$ by Riemann-Lebesgue's Lemma, and the substitution
$u=x\La$ in the second integral leads to
\begin{eqnarray*}
\lefteqn{\int_0^s
x^{-(c-1)}\varphi_1(x)\,dx}\\
&=&\La^{c-2}\int_0^{s\La}u^{-(c-1)}\varphi_1(u/\La)\,du\\
&=&-\La^{c-1}\int_0^{s\La} u^{-(c-1)}e^{-u}\,du -\La^{c-1}\int_0^{s\La} u^{-(c-1)}e^{-u}\left(\Gamma(1-u/\La)^c-1\right)\,du\\
&-&
c\La^{c-2}\int_0^{s\La}u^{-(c-1)}e^{-u}\Gamma(1-u/\La)^c\Psi(1-u/\La)\,du\\
&=&-\La^{c-1}\Gamma(2-c)+\mathcal O(\La^{c-2}).
\end{eqnarray*}
Using that
$$
 \left(-\frac{\sin(\pi(c-1))}{(c-1)\pi}\right)\left(-\La^{c-1}\Gamma(2-c)\right)=\frac{\La^{c-1}}{\Gamma(c)}
$$
by Euler's reflection formula, we see that \eqref{eq:asat0} holds.

\bigskip
{\bf Case 3.} We now assume $c>2$.

We perform the change of variable $w=\La z$ in \eqref{eq:cont1} and obtain
$$
e_c(t)=\frac{\La^{c-1}}{2\pi i}\int_{V_{-\La}}(-w)^{-c} e^{-w}\Gamma(1-w/\La)^c\,dw.
$$
  Using Cauchy's integral theorem, we can shift the contour $V_{-\La}$
  to $V_{-1}$  as the integrand is holomorphic in the vertical strip between both paths and exponentially small at both extremes of that vertical strip. Then,
$$
e_c(t)=\frac{\Lambda^{c-1}}{2\pi i}\int_{V_{-1}}(-w)^{-c} e^{-w}\Gamma\left(1-w/\La\right)^c\,dw.
$$
For any holomorphic function $h$ in a domain $G$ which is  star-shaped
with respect to $0$ we have
$$
h(z)=h(0)+z\int_0^1 h'(uz)\,du,\quad z\in G.
$$ 
If this is applied to $G=\mathbb C\setminus[1,\infty)$ and
$h(z)=\Gamma(1-z)^c$ we find
\begin{equation}\label{eq:taylor}
\Gamma(1-w/\La)^c=1-\frac{cw}{\La}\int_0^1 \Gamma(1-uw/\La)^c\Psi(1-uw/\La)\,du.
\end{equation}
 Defining 
$$
R(w)=\int_0^1 \Gamma(1-uw/\La)^c\Psi(1-uw/\La)\,du,
$$
we get
\begin{equation}\label{eq:taylor1}
e_c(t)={\Lambda^{c-1}}{2\pi i}\int_{V_{-1}} (-w)^{-c} e^{-w}dw +
\frac{c\Lambda^{c-2}}{2\pi i}\int_{V_{-1}}(-w)^{1-c} e^{-w}R(w)dw.
\end{equation}
For any $w\in V_{-1}, 0\le u\le 1$ and for $\La\ge 1$ we have that
 $1-uw/\La\in \Omega$, where $\Omega$ is the closed vertical strip
 located between the vertical lines $V_1$ and $V_2$. Because
$\Gamma(z)^c\Psi(z)$ is  continuous in $\Omega$ and exponentially small at the upper and lower limits of $\Omega$, the function
$R(w)$ is bounded for $w\in V_{-1}$ by a constant independent of $\La\ge 1$. Therefore,
$$
\frac{c\Lambda^{c-2}}{2\pi i}\int_{V_{-1}} (-w)^{1-c} e^{-w}R(w)dw=\mathcal O(\Lambda^{c-2}),
$$
where we use that $(-w)^{1-c}e^{-w}$ is integrable over $V_{-1}$
because $c>2$.

On the other hand, in the first integral of \eqref{eq:taylor1}, the contour $V_{-1}$ may be deformed to a Hankel contour 
$$
\mathcal H:=\{x-i \mid x=\infty\ldots  0\}\cup  \{e^{i\theta}\mid \theta=-\pi/2 \ldots -3\pi/2\}
\cup\{x+i\mid x=0\ldots \infty\}
$$
surrounding $[0,\infty)$, and the integral over $\mathcal H$ is Hankel's integral representation of the inverse of the Gamma function:
$$
\frac{1}{2\pi i}\int_{\mathcal H} (-w)^{-c} e^{-w}dw=\frac{1}{\Gamma(c)}.
$$
Therefore,  when we join everything, we obtain that for $c>2$:
$$
e_c(t)=\frac{(\log(1/t))^{c-1}}{\Gamma(c)}+\mathcal O((\log(1/t))^{c-2}), \quad t\to 0.
$$

\bigskip
{\bf Case 4.} $c=1,c=2$.

These cases are easy since $e_1(t)=e^{-t}$ and
$e_2(t)=2K_0(2\sqrt{t})$. $\quad\square$

\begin{rem} {\rm The behaviour of $e_c(t)$ for $t\to 0$ can be
    obtained from \eqref{eq:fund1} using the residue theorem when $c$
    is a natural number. In fact, in this case $\Gamma(-z)^c$ has a
    pole of order $c$ at $z=0$, and a shift of the contour $V_{-1}$ to
    $V_s$, where $0<s<1$, has to be compensated by a residue, which
    will give the behaviour for $t\to 0$.}
\end{rem}

\medskip
{\bf Acknowledgment:} The authors want to thank Nico Temme for his indications about the asymptotics of the integral \eqref{eq:U1}.

\medskip
C. Berg, Department of Mathematical Sciences, University of
Copenhagen, Universitetsparken 5,
2100 Copenhagen {\O}, Denmark

email: berg@math.ku.dk

\medskip
J. L. L{\'o}pez, Departamento de Ingener\'{\i}a Matem{\'a}tica e
Inform{\'a}tica, Universidad P{\'u}blica de Navarra, 31006 Pamplona,
Spain

email: jl.lopez@unavarra.es
 
\end{document}